\documentclass[12pt,a4wide]{article}
\usepackage[leqno]{amsmath}
\usepackage[all]{xy}
\usepackage[mathscr]{euscript}
\usepackage{amssymb}
\usepackage{amscd}
\usepackage[pdftex]{color,graphicx}

\newtheorem{lemma}{Lemma}[section]
\newtheorem{theorem}[lemma]{Theorem}
\newtheorem{prop}[lemma]{Proposition}
\newtheorem{claim}[lemma]{Claim}
\newtheorem{cor}[lemma]{Corollary}
\newtheorem{remark}[lemma]{Remark}

\newtheorem{conj}[lemma]{Conjecture}
\newtheorem{question}[lemma]{Question}

\newcommand{\pf}{\noindent{\em Proof: }}
\newcommand{\epf}{\hfill\hbox{\rule{3pt}{6pt}}\\}

\newcommand{\forme}[1]{}

\begin{document}

\title{On the limit points of the smallest eigenvalues of regular graphs}

\author{ Hyonju Yu\\
{\small {\tt lojs4110@postech.ac.kr}}\\
{\footnotesize{Department of Mathematics,  POSTECH, Pohang 790-785, South Korea}}}

\date{\today}

\maketitle
\begin{abstract}
In this paper, we give infinitely many examples of (non-isomorphic) connected $k$-regular graphs with smallest eigenvalue in half open interval $[-1-\sqrt2, -2)$ and also infinitely many examples of (non-isomorphic) connected $k$-regular graphs with smallest eigenvalue in half open interval $[\alpha_1, -1-\sqrt2)$ where $\alpha_1$ is the smallest root$(\approx -2.4812)$ of the polynomial $x^3+2x^2-2x-2$. From these results, we determine the largest and second largest limit points of smallest eigenvalues of regular graphs less than $-2$. Moreover we determine the supremum of the smallest eigenvalue among all connected 3-regular graphs with smallest eigenvalue less than $-2$ and we give the unique graph with this supremum value as its smallest eigenvalue.
\end{abstract}

\section{Introduction}

Line graphs are examples of graphs with smallest eigenvalue at least $-2$ (where the eigenvalues of a graph are the eigenvalues of its adjacency matrix). In 1976, P. J. Cameron et al. characterized the graphs with smallest eigenvalue at least $-2$ \cite{180}. They showed (using root lattices):

\begin{theorem}(Cf. \cite[Theorem 3.12.1]{bcn})
A connected graph $G$ with smallest eigenvalue at least $-2$ is a generalized line graph (for definition, see next section) or the number of vertices of $G$ is at most 36.
\end{theorem}

A. J. Hoffman studied graphs with smallest eigenvalue less than $-2$. He showed that the largest limit point smaller than $-2$ for the smallest eigenvalues of graphs is $-1-\sqrt 2$ \cite{Hoff}. More precisely, he showed:

\begin{theorem}\label{Hoff}
Let $\widehat{\theta}_k$ be the supremum of the smallest eigenvalues of graphs with minimal valency $k$ and smallest eigenvalue $<-2$. Then $(\widehat{\theta}_k)_k$ forms a monotone decreasing sequence with limit $-1-\sqrt 2$.
\end{theorem}

A consequence of the results of R. Woo and A. Neumaier is that $\alpha_1$ (where $\alpha_1$ is the smallest root$(\approx -2.4812)$ of the polynomial $x^3+2x^2-2x-2$) is the largest limit point of the smallest eigenvalues of graphs with smallest eigenvalues less than $-1-\sqrt 2$ \cite{woo}. More precisely:

\begin{theorem}\label{woolim}
Let $\widehat{\sigma}_k$ be the supremum of the smallest eigenvalues of graphs with minimal valency $k$ and smallest eigenvalue $<-1-\sqrt 2$. Then $(\widehat{\sigma}_k)_k$ forms a monotone decreasing sequence with limit $\alpha_1$.
\end{theorem}

We will give proofs of Theorems \ref{Hoff} and \ref{woolim} in Section 2.

In this paper, we study the limit points of the smallest eigenvalues of regular graphs. Clearly $-2$ is a limit point for regular graphs. Also $\frac{-3-\sqrt 5}{2}$ is a limit point for regular graphs as the cartesian product of a pentagon and a complete graph with at least 2 vertices has smallest eigenvalues $\frac{-3-\sqrt 5}{2}$. But it is not known whether this is the largest limit point smaller then $-2$ for regular graphs.

In this paper, we will construct $k$-regular graphs which have the smallest eigenvalues in $[-1-\sqrt2, -2)$ for all $k$. And as a corollary, we will show that for regular graphs, the largest limit point of regular graphs with smallest eigenvalue less than $-2$ is also $-1-\sqrt 2$. Our results are:

\begin{theorem}\label{infimany}
Let $k$ be a positive integer at least 3. Then there are infinitely many $k$-regular graphs which have smallest eigenvalue in $[-1-\sqrt2, -2)$. Also, there are infinitely many $3$-regular triangle free graphs which have the smallest eigenvalue in $[-1-\sqrt2, -2)$.
\end{theorem}

\begin{cor}\label{short}
Let $\hat{\eta}_k$ be the supremum of the smallest eigenvalues of $k$-regular graphs with smallest eigenvalue $<-2$. Then $$\lim_{k\rightarrow\infty}(\hat{\eta}_k)_k=-1-\sqrt{2}.$$
\end{cor}

Furthermore, for enough large $k$, we will construct $k$-regular graphs which have smallest eigenvalues in $[\alpha_1, -1-\sqrt2)$ where $\alpha_1$ is the smallest root$(\approx -2.4812)$ of the polynomial $x^3+2x^2-2x-2$. And as a corollary, we will show that the largest limit point of regular graphs with smallest eigenvalue less than $-1-\sqrt 2$ is also $\alpha_1$, as shown for general graphs by R. Woo and A. Neumaier.

\begin{theorem}\label{infimanyhh}
There exist $N$ such that there are infinitely many $k$-regular graphs which have smallest eigenvalue in $[\alpha_1, -1-\sqrt2)$ for $k\geq N$.
\end{theorem}

\begin{cor}\label{short1}
Let $\hat{\xi}_k$ be the supremum of the smallest eigenvalues of $k$-regular graphs with smallest eigenvalue $<-1-\sqrt 2$. Then $$\lim_{k\rightarrow\infty}(\hat{\xi}_k)_k=\alpha_1$$
where $\alpha_1$ is the smallest root$(\approx -2.4812)$ of the polynomial $x^3+2x^2-2x-2$.
\end{cor}

Theorems \ref{infimany} disproves the following conjecture by R. Woo and A. Neumaier.

\begin{conj}\cite[p.590, Conjecture (ii)]{woo}\label{conj}
A regular graph with smallest eigenvalue $\geq -1-\sqrt 2$ and sufficiently large valency is a line graph or a cocktail party graph.
\end{conj}

In order to show Theorems \ref{infimany} and \ref{infimanyhh}, we will use Hoffman graphs, which were introduced by R. Woo and A. Neumaier\cite{woo}.

In the next section we will give definitions and preliminaries. Also in this section we will recall Hoffman graphs and some of their basic theory. In Section 3, we will show Theorems \ref{infimany}, \ref{infimanyhh} and Corollary \ref{short}, \ref{short1}. In Section 4, we will determine $\hat{\eta}_3$. In the last section, we will give concluding remarks.

\section{Definitions and preliminaries}

All the graphs considered in this paper are finite, undirected and without multiple edges.
A graph $G$ is an ordered pair $(V(G), E(G))$ of vertex set $V=V(G)$ and edge set $E=E(G)$, where $E(G)$ consists of unordered pairs of two adjacent vertices, so we can consider $E(G)\subseteq\binom{V(G)}{2}$.

For two graph $G$ and $H$ such that $V(H)\subseteq V(G)$ and $\{x,y\}\in E(H)$ if and only if $\{x,y\}\in E(G)$ for all $x, y \in V(H)$, we call that $H$ is induced subgraph of $G$ and $G$ is super graph of $H$.

The adjacency matrix $A(G)$ of a graph $G$ is the $(0,1)$-matrix whose rows and columns are indexed by the vertex set $V(G)$ and the $(x,y)$-entry is $1$ whenever $x$ is adjacent to $y$ and $0$ otherwise. The eigenvalues of $G$ are the eigenvalues of its adjacency matrix. A graph $G$ is $k$-regular if for any vertex $x$, the number of edges which contain $x$, we call valency of $x$, is $k$. A graph $G$ is regular if it is $k$-regular for some integer $k$. A graph $G$ is $(k_1, k_2)$-semiregular if there is a partition $\{A\neq\emptyset, B\neq\emptyset\}$ of vertex set such that every vertex in $A$ has valency $k_1$ and every vertex in $B$ has valency $k_2$. We denote $\delta(G)$ the minimum valency of $G$.

For a given graph $G$, the line graph of $G$, $\L(G)$, has vertex set $E(G)$ and $\{x, y\}$ is adjacent to edge $\{v, w\}$ if and only if $|\{x, y\}\cap \{v, w\}|=1$. A graph $G$ is called a line graph if $G=L(H)$ for some graph $H$.

For a line graph $G$, say $G=L(H)$, we obtain the relation $C^T C=A(G)+2I$ where $C$ is a vertex-edge incidence matrix of $H$ and $I$ is an identity matrix. This relation implies that the smallest eigenvalue of a line graph is at least $-2$.

Line graphs can be recognized combinatorially as J. Krausz \cite{kra} showed:

\begin{theorem}(Cf. \cite[Theorem 2.1.1]{line})\label{linechar}
A graph is a line graph if and only if its edge set can be partitioned into non-trivial cliques such that:
\begin{itemize}
\item[(i)] two cliques have at most one vertex in common;
\item[(ii)] each vertex is in at most two cliques.
\end{itemize}
\end{theorem}

In this paper, we are mainly interested in regular graphs. P. J. Cameron et al. \cite{180} characterized the connected regular graphs with smallest eigenvalue at least $-2$. Their result:

\begin{theorem}(Cf. \cite[Theorem 3.12.2]{bcn})\label{cla1}
Let $G$ be a connected regular graph with $\nu$ points, valency $k$, and smallest eigenvalue at least $-2$. Then one of the following holds:
\begin{itemize}
\item[(i)] $G$ is the line graph of a regular or a bipartite semiregular connected graph.
\item[(ii)] $\nu=2(k+2)\leq 28$ and $G$ is an induced subgraph of $E_7(1)$.
\item[(iii)] $\nu=\frac{3}{2}(k+2)\leq 27$ and $G$ is an induced subgraph of Schl\"{a}fli graph.
\item[(iv)] $\nu=\frac{4}{3}(k+2)\leq 16$ and $G$ is an induced subgraph of Clebsch graph.
\item[(v)] $\nu=k+2$ and $G\cong K_{m\times 2}$ for some $m\geq3$.
\end{itemize}
\end{theorem}

This theorem means that there are only finitely many connected regular graphs with smallest eigenvalue at least $-2$ which are neither line graphs nor Cocktail party graphs.

\subsection{Hoffman graphs}

Hoffman graphs can be used as a tool to construct graphs with smallest eigenvalue at least a fixed number. In this subsection, we will define Hoffman graphs (as introduced by R. Woo and A. Neumaier) and discuss their basic theory \cite{woo}.

A Hoffman graph $\mathfrak{H}$ is a pair $(H=(V, E), \ \mu : V\rightarrow \{f, s\})$ satisfying the following conditions:
\begin{itemize}
\item[(i)] H is a graph;
\item[(ii)] every vertex with label $f$ is adjacent to at least one vertex with label $s$;
\item[(iii)] vertices with label $f$ are pairwise non-adjacent.
\end{itemize}

We call a vertex with label $s$ a slim vertex, and a vertex with label $f$ a fat vertex. We denote by $V_s=V_s(\mathfrak H)$ (resp. $V_f=V_f(\mathfrak H)$) the set of slim (resp. fat) vertices of $\mathfrak H$. The subgraph of a Hoffman graph $\mathfrak H$ induced on $V_s(\mathfrak H)$ is called the slim graph of $\mathfrak H$.

Let $\mathfrak{H}$ be a Hoffman graph and let $\mathfrak{H}^1$ and $\mathfrak{H}^2$ be two non-empty induced Hoffman subgraphs of $\mathfrak{H}$. The  Hoffman graph $\mathfrak{H}$ is said to be the sum of $\mathfrak{H}^1$ and $\mathfrak{H}^2$, written as $\mathfrak{H}=\mathfrak{H}^1\oplus\mathfrak{H}^2$, if the following conditions are satisfied:
\begin{itemize}
\item[(i)] $V(\mathfrak{H})=V(\mathfrak{H}^1)\cup V(\mathfrak{H}^2)$;
\item[(ii)] $\{V_s(\mathfrak{H}^1), V_s(\mathfrak{H}^2)\}$ is a partition of $V_s(\mathfrak{H})$;
\item[(iii)] if $x\in V_s(\mathfrak{H}^1)$ and $y\in V_s(\mathfrak{H}^2)$, then $x$ and $y$ have at most one common fat neighbor and they have one if and only if they are adjacent;
\item[(iv)] if $x\in V_s(\mathfrak{H}^i)$ and $F\in V_{f}(\mathfrak{H})$ are adjacent in $\mathfrak{H}$, then $F\in V_{f}(\mathfrak{H}^i)$ and $x$ and $F$ are adjacent in $\mathfrak{H}^i$ for $i=1, 2$.
\end{itemize}

It easily follows that the sum defined above is associative, in the sense that if $\mathfrak{H}=\mathfrak{H}^1\oplus (\mathfrak{H}^2\oplus \mathfrak{H}^3)$, then $\mathfrak{H}=(\mathfrak{H}^1\oplus \mathfrak{H}^2)\oplus \mathfrak{H}^3$ and vice versa. Instead of $\mathfrak{H}=((\ldots(\mathfrak{H}^1\oplus \mathfrak{H}^2)\oplus \mathfrak{H}^3) \ldots \mathfrak{H}^n)$, we write $\mathfrak{H}=\bigoplus_{i=1}^n \mathfrak{H}^i$.

Let $\mathfrak{A}$ be a family of Hoffman graphs. An $\mathfrak{A}$-line graph is an induced subgraph of the slim graph of a Hoffman graph $\mathfrak{H}=\bigoplus_{i=1}^s \mathfrak{H}^i$ for some $s$.

In particular, it can be easily checked that a line graph is a $\{\mathfrak{H_2}\}$-line graph where $\mathfrak{H_2}$ is as shown in {\em Figure 1} at the end of this subsection. We define a generalized line graph as a $\{\mathfrak{H_2}, \mathfrak{H_3}\}$-line graph where $\mathfrak{H_3}$ is also shown in {\em Figure 1}. This definition is equivalent with original definition in \cite[p. 106]{bcn}.

Now, we define eigenvalues of Hoffman graphs.
For a Hoffman graph $\mathfrak H=(H, \mu)$, let $A$ be the adjacency matrix of $H$,
$$A=\left[
\begin{array}{cc}
 A_s  & C \\
 C^T & 0
\end{array} \right]$$
where $A_s$ is the adjacency matrix of the slim graph of $\mathfrak{H}$. Eigenvalues of $\mathfrak H$ are defined as the eigenvalues of the real symmetric matrix $B=B(\mathfrak{H}):= A_s - C C^T$. By $\lambda_{min}(\mathfrak H)$, we denote the smallest eigenvalue of $\mathfrak H$.

From above definitions, we can regard any ordinary graph $H$ as a Hoffman graph $\mathfrak H=(H, \mu)$ satisfing $\mu(x)=s$ for all $x\in V(H)$, and then the eigenvalues of $H$ are exactly the same as the eigenvalues of $\mathfrak{H}$.

We denote by $\lambda_{min}=\lambda_{min}(G)$ the smallest eigenvalue of a given graph $G$. {\em Figure 1} shows that the smallest eigenvalues of some Hoffman graphs where $\alpha_0 = -1-\sqrt2$ and $\alpha_1\approx -2.4812$ is the smallest root of $x^3 +2x^2-2x-2$.

The following theorem shows the relation between the smallest eigenvalue of a Hoffman graph and any of its induced subgraph, in particular, its slim graph.

\begin{theorem}(Cf. \cite[Theorem 3.2 and 3.7]{woo})\label{sinter}
\begin{itemize}
\item[(i)] If $\mathfrak G$ is an induced subgraph of a Hoffman graph $\mathfrak H$, then $\lambda_{min}(\mathfrak{G})\geq \lambda_{min}(\mathfrak{H})$ holds.
\item[(ii)] If $\mathfrak{H}=\mathfrak{H_1}\oplus\mathfrak{H_2}$, then $\lambda_{min}(\mathfrak{H})=\min\{\lambda_{min}(\mathfrak{H_1}), \lambda_{min}(\mathfrak{H_2})\}$.
\end{itemize}
\end{theorem}

To consider the limit points of smallest eigenvalues, we discuss the following proposition which was shown by A. J. Hoffman and A. M. Ostrowski. (For a proof, see \cite{JKMT}).

\begin{prop}(Cf. \cite[Proposition 5.3]{woo})\label{Hoffthm}
Let $\mathfrak{H}$ be a Hoffman graph and denote by $\mathfrak{H}^{(n)}$ the clique extension obtained by replacing all fat vertices by slim $n$-cliques, all of whose vertices are adjacent to the neighbors of the corresponding fat vertices. Then
$$\lim_{n\rightarrow \infty} \lambda_{min}(\mathfrak{H}^{(n)})=\lambda_{min}(\mathfrak{H}).$$
\end{prop}

The following theorem gives a structure theory for graph with smallest eigenvalue a little bit smaller than $-2$.

\begin{theorem}\label{Hoff1wootheorem}
\begin{itemize}
\item[(i)] (Cf. \cite[Theorem 1.1]{Hoff})
For $-1-\sqrt2 <\lambda \leq -2$, there exists a number $f(\lambda)$ such that if $G$ is a connected graph with minimum valency of $G$ is at least $f(\lambda)$ and $\lambda_{min}(G)\geq \lambda$, then $\lambda_{min}(G)=-2$ and $G$ is a generalized  line graph.
\item[(ii)] (Cf. \cite[Theorem 5.1]{woo})
Let $\lambda\leq -1$ be a real number larger than $\alpha_1\approx-2.4812$, the smallest root of the polynomial $x^3+2x^2-2x-2$. Then there exists a number $\kappa_{\lambda}$ such that every graph with minimum valency at least $\kappa_{\lambda}$ such that its smallest eigenvalue $\lambda_{min}$ at least $\lambda$ is an $\{{\mathfrak H}_2, {\mathfrak H}_5, {\mathfrak H}_7, {\mathfrak H}_9\}$-line graph. In particular, $\lambda_{min}\geq -1-\sqrt 2$.
\end{itemize}
\end{theorem}

{\bf Proofs of Theorem \ref{Hoff} and Theorem \ref{woolim}} \\

Proof of Theorem \ref{Hoff}:
By Theorem \ref{Hoff1wootheorem}(i), we have $\lim_{k\rightarrow \infty}\widehat{\theta}_k\leq -1-\sqrt2$. By Proposition \ref{Hoffthm} and the Hoffman graph $\mathfrak{H}_{9}$ in {\em Figure 1}, we construct a sequence of graphs $(G_n)_{n=2}^{\infty}$ such that $(\lambda_{min}(G_n))_n$ is a monotone decreasing sequence and $\lim_{n\rightarrow\infty}\delta(G_n)=\infty$. This shows Theorem \ref{Hoff}. \\

Proof of Theorem \ref{woolim}:
The proof is similar except that we need to change Theorem \ref{Hoff1wootheorem}(i) to Theorem \ref{Hoff1wootheorem}(ii), $\mathfrak{H}_{9}$ to $\mathfrak{H}_{WN}$ and $-1-\sqrt2$ to $\alpha_1$. This shows Theorem \ref{woolim}.
\epf

\begin{center}
\begin{figure}[h]
\includegraphics[width=0.8\textwidth]{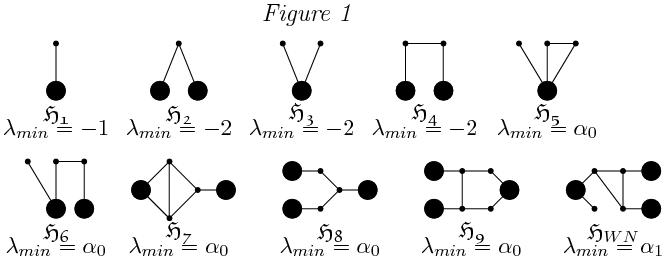}
\end{figure}
\end{center}

\section{Proofs of the main theorems}

In this section, we will give the proofs for Theorems \ref{infimany}, \ref{infimanyhh} and Corollaries \ref{short}, \ref{short1}.

\subsection{Proof of the Theorem \ref{infimany}}\label{3}

Let $k$ be a positive integer at least 3.

\begin{remark}\label{semi}
There exist a connected $(k-1, k)$-semi-regular graph with $a(2k-1)$ vertices for all positive integer $a$.
\end{remark}

Let $B_{k, k-1}$ be a semi-regular bipartite graph of partition $R$ with valency $k$ and $Y$ with valency $k-1$. For the graph $B_{k, k-1}$, we are going to construct connected $k$-regular graph $G(B_{k, k-1})=:G_k$ with $\lambda_{min}(G_k))\geq -1-\sqrt2$. The graph $G_k$ will be the slim graph of the Hoffman graph $\mathfrak{G}(B_{k, k-1})=:\mathfrak{G}_k$ as defined below. Let

\begin{center}
$V_s(\mathfrak{G}_k)=\{x_{t}^{(i,j)_k}\mid \{i, j\}\in E(B_{k, k-1}), i\in R , j\in Y$ and $t=1, 2, 3\}$  \\
and $V_f(\mathfrak{G}_k)=\{F_i, E_i\mid i\in R\}\cup \{D_j\mid j\in Y\}$. \\
\end{center}

So, $\sharp V_s(\mathfrak{G}_k)=3k(\sharp R)=3(k-1)(\sharp Y)$ and $\sharp V_f(\mathfrak{G}_k)=2\sharp R +\sharp Y$. Define

\begin{center}
\begin{figure}[h]
\includegraphics[width=0.8\textwidth]{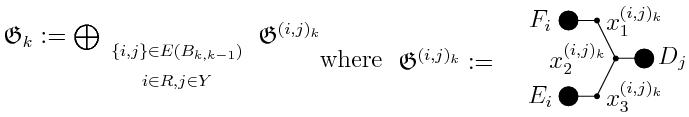}
\end{figure}
\end{center}

and let $G_k$ be the slim graph of $\mathfrak{G}_k$.

Then the following claim holds:
\begin{claim}
\begin{itemize}
\item[(i)] $G_k$ is a $k$-regular graph;
\item[(ii)] $\lambda_{min}(G_k)\geq \lambda_{min}(\mathfrak{G}_k)=\lambda_{min}(\mathfrak{H}_{8})=-1-\sqrt2$ for all $\{i, j\}\in E(B_{k, k-1})$ where $i\in R$ and $j\in Y$;
\item[(iii)] $G_k$ is neither a line graph or a Cocktail party graph.
\end{itemize}
\end{claim}

\pf
$(i)$ : By construction of $G_k$. \\
$(ii)$ : By Theorem \ref{sinter}. \\
$(iii)$ : For any two edges $\{i_1, j\}$ and $\{i_2, j\}$ of $B_{k, k-1}$ with common vertex $j$ in $Y$, the induced subgraph on $\{x_1^{(i_1,j)_k}, x_2^{(i_1,j)_k}, x_3^{(i_1,j)_k}, x_2^{(i_2,j)_k}\}$ is 3-claw. By Theorem \ref{linechar}, this implies $(iii)$.

\epf

By Theorem \ref{cla1} and Remark \ref{semi}, we obtain infinitely many $k$-regular graphs which have smallest eigenvalue in the half open interval $[-1-\sqrt2 , -2)$. \\

Now, we consider 3-regular triangle free graphs.

Let $C_{2n}$ be a $2n$-gon with vertex set $\{x_{2n+1}=x_1, x_2, \ldots , x_{2n}\}$ and edge set $\{\{x_i, x_{i+1}\}\mid i=1, 2, \ldots , 2n\}$. For the graph $C_{2n}$, we are going to construct $3$-regular triangle free graph $G(C_{2n})=:G'$ with $\lambda_{min}(G')\geq -1-\sqrt2$. The graph $G'$ will be the slim graph of the Hoffman graph $\mathfrak{G}(C_{2n})=:\mathfrak{G}'$. Let

\begin{center}
$V_s(\mathfrak{G}')=\{x_{t}^{(2i-1,2i)}, x_{t}^{(2i-1,2i-2)}\mid i=1, 2, \ldots , n$ and $t=1, 2, 3,4\}$ \\
and $V_f(\mathfrak{G}')=\{F_{2i-1}, E_{2i-1}\mid i=1, 2, \ldots , n\}\cup \{D_{2i}\mid i=1, 2, \ldots , n\}$.
\end{center}

Define $\mathfrak{G}'=\bigoplus_{i\in\{1, 2, \ldots , n\}} \left(\mathfrak{G}^{(2i-1, 2i)}\bigoplus \mathfrak{G}^{(2i-1, 2i-2)}\right)$ where

\vspace{-0.5cm}

\begin{center}
\begin{figure}[h]
\includegraphics[width=0.8\textwidth]{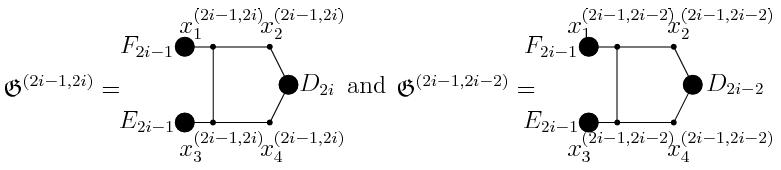}
\end{figure}
\end{center}

\vspace{-0.5cm}

and let $G'$ be the slim graph of $\mathfrak{G}'$.

Then the following claim holds:
\begin{claim}
\begin{itemize}
\item[(i)] $G'$ is a $3$-regular graph;
\item[(ii)] $\lambda_{min}(G')\geq \lambda_{min}(\mathfrak{G}')=\lambda_{min}(\mathfrak{H}_{9})=-1-\sqrt2$;
\item[(iii)] Neighbors of a fixed vertex are not adjacent each other in $G'$;
\item[(iv)] $G'$ is neither a line graph nor a Cocktail party graph;
\item[(v)] $G'$ is triangle free.
\end{itemize}
\end{claim}

\pf
$(i)$ : By construction of $G'$. \\
$(ii)$ : By Theorem \ref{sinter}. \\
$(iii)$ : By symmetry, we only need to consider the neighbors of $x_1^{(1,2)}$ and $x_2^{(1,2)}$, those are $\{x_2^{(1,2)}, x_3^{(1,2)}, x_1^{(1,2n)}\}$ and $\{x_1^{(1,2)}, x_2^{(3,2)}, x_4^{(3,2)}\}$. \\
And $(iii)$ implies $(iv)$ and $(v)$ .
\epf

By Theorem \ref{cla1}, we obtain infinitely many $3$-regular triangle free graphs which have smallest eigenvalue in the half open interval $[-1-\sqrt2 , -2)$.

\epf

\subsection{Proof of Theorem \ref{infimanyhh}}\label{2}

By Proposition \ref{Hoffthm}, we can take a positive integer $N$ such that $\lambda_{min}(\mathfrak{H}_{WN}^{(n-3)})< -1-\sqrt2$ for all $n\geq N$. \\

Let $k\geq N$ and $a$ be positive integers. \\

We are going to construct a $k$-regular graph $G_k$ such that $\alpha_1 \leq\lambda_{min}(G_k)< -1-\sqrt 2$. The graph $G_k$ will be the slim graph of a Hoffman graph $\mathfrak{G}_k$ with $2a(2k^2-4k+1)$ fat vertices and $4ak(k-1)(k-2)$ slim vertices. \\

Let $P(k)$, $Q(k)$ and $R(k)$ be three partitions of $\{1, 2, \ldots , ak(k-1)(k-2)\}=:[ak(k-1)(k-2)]$ such that
\begin{center}
${P(k)}=\{{p(k)}_i \mid i=1, 2, \ldots , ak(k-1) \}\subseteq \binom{[ak(k-1)(k-2)]}{k-2}$ \\
${Q(k)}=\{{q(k)}_j \mid j=1, 2, \ldots , a(k-1)(k-2)\}\subseteq \binom{[ak(k-1)(k-2)]}{k}$ \\
${R(k)}=\{{r(k)}_l \mid l=1, 2, \ldots , ak(k-2)\}\subseteq \binom{[ak(k-1)(k-2)]}{k-1}$.
\end{center}

And let

\begin{center}
${I(k)}=\{(i,j,l)_m \mid m\in {p(k)}_i\cap {q(k)}_j\cap {r(k)}_l$ for $1\leq m\leq ak(k-1)(k-2)\}$.
\end{center}

Now, we define Hoffman graph $\mathfrak{G}_k$ with fat vertex set

\begin{center}
$V_f(\mathfrak{G}_k)=\{F_i \mid i=1, \ldots , ak(k-1) \}$ \\ $\cup \{ E_j \mid j=1, \ldots , a(k-1)(k-2)\} \cup\{ D_l, C_l \mid l=1, \ldots , ak(k-2)\}$
\end{center}

and slim vertex set

\begin{center}
$V_s(\mathfrak{G}_k)=\{x_{t}^{(i,j,l)_m}\mid (i,j,l)_m\in I(k), \ t=1, 2, 3, 4\}$
\end{center}

as following:

\vspace{-0.5cm}

\begin{center}
\begin{figure}[h]
\includegraphics[width=0.8\textwidth]{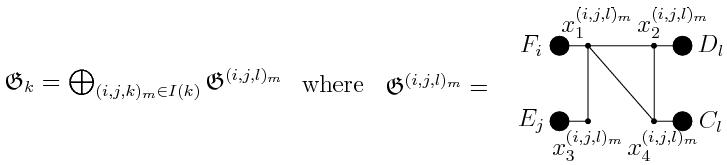}
\end{figure}
\end{center}

\vspace{-0.5cm}

and let $G_k$ be the slim graph of $\mathfrak{G}_k$.

For $p(k)_i=\{\alpha(1), \alpha(2), \ldots , \alpha(k-2)\}$, there exist $j(s)$ and $l(s)$ such that $\alpha(s)\in q(k)_{j(s)}\cap r(k)_{l(s)}$. This means that $(i, j(s), l(s))_{\alpha(s)}\in I(k)$ for $s=1, 2, \ldots , k-2$ and there are $k-2$ slim neighbors of fat vertex $F_i$. By a similar argument, there are $k$, $k-1$ and $k-1$ slim neighbors of fat vertex $E_j$, $D_l$ and $C_l$, respectively. Moreover we have the following claim:

\begin{claim}
\begin{itemize}
\item[(i)] $G_k$ is a $k$-regular graph;
\item[(ii)] $\lambda_{min}(G_k)\geq \lambda_{min}(\mathfrak{G}_k)=\lambda_{min}(\mathfrak{H}_{WN})=\alpha_1$;
\item[(iii)] $\lambda_{min}(\mathfrak{G}_k)< -1-\sqrt2$.
\end{itemize}
\end{claim}

\pf
$(i)$ : By construction of $G_k$. \\
$(ii)$ : By Theorem \ref{sinter}. \\
$(iii)$ : Consider the induced subgraph on $x_1^{(i,j,k)_m}, x_2^{(i,j,k)_m}, x_3^{(i,j,k)_m}, x_4^{(i,j,k)_m}$ and its neighbors for some fixed $(i,j,k)_m \in I(k)$. Then this induced subgraph contains $\mathfrak{H}_{WN}^{(k-3)}$. So, $\lambda_{min}(G_k)\leq \lambda_{min}(\mathfrak{H}_{WN}^{(k-3)})<-1-\sqrt2$ as $k\geq N$.

\epf

It is not so difficult to see that without loss of generality we may assume that $G_k$ is connected for any positive integer $a$. We have constructed connected $k$-regular graphs on $4ak(k-1)(k-2)$ vertices which have the smallest eigenvalue in the half open interval $[\alpha , -1-\sqrt2)$ for $k\geq N$. Therefore, our result is proved.

\epf

\subsection{Proofs of Corollary \ref{short} and \ref{short1}}\label{3}

By definition of $\hat{\theta_k}$ and $\hat{\eta_k}$, $\hat{\theta_k}\geq \hat{\eta_k}$. This implies that $\lim_{k\rightarrow\infty}\hat{\eta_k}\leq -1-\sqrt2$ as Theorem \ref{Hoff}. Since Theorem \ref{infimany} means that $\lim_{k\rightarrow\infty}\hat{\eta_k}\geq -1-\sqrt2$, we obtain Corollary \ref{short}. As similar argument, Theorem \ref{woolim} and Theorem \ref{infimanyhh} imply Corollary \ref{short1}.

\section{The value $\hat{\eta}_3$}

In order to characterize $\hat{\eta}_3$, we use following lemma which follows immediately from Theorem \ref{linechar}.

\begin{lemma}\label{line}
(i) A graph with a 3-claw as an induced subgraph is not a line graph.  \\
(ii) Let $G$ be a 3-regular graph which is not $K_4$. Then $G$ is a line graph if and only if $N(x)\cong K_1 \cup K_2$ for all $x$ in $V(G)$.
\end{lemma}

The following theorem determines $\hat{\eta}_3$ and gives the graph such that the smallest eigenvalue is $\hat{\eta}_3$.

\begin{theorem}\label{limit}
Let $\beta\approx-2.0391$ be the smallest root of $x^6-3x^5-7x^4+21x^3+13x^2-35x-4$. Then $\hat{\eta}_3=\beta$ and

\vspace{-0.5cm}

\begin{center}
\begin{figure}[h]
~~~~~~~~~~~~~~~~~~~~~~~~~~~~~~~~~~~~~\includegraphics[width=0.2\textwidth]{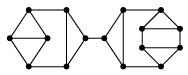}
\end{figure}
\end{center}

\vspace{-1cm}

is the unique connected 3-regular graph with smallest eigenvalue $\beta$.
\end{theorem}

\pf
As the graph, shown in the theorem, has smallest eigenvalue $\beta\approx-2.0391$, we only need to show that the only connected $3$-regular graph with $\beta\leq\lambda_{min}<-2$ is the graph of this theorem. Let $G$ be a connected $3$-regular graph with smallest eigenvalue $\beta\leq \lambda_{min}(G)<-2$. Then $G$ contains a $3$-claw or a $K_{2,1,1}$ as induced subgraph by Lemma \ref{line}. We start with the assumption that $G$ contains a $K_{2,1,1}$, that is the graph $A1$ of Appendix 1. Now we proceed by adding vertices to a vertex that does not have degree $3$ and consider all possible cases . If such a graph has smallest eigenvalue smaller than $\beta$, then we know that it can not be an induced subgraph of $G$ by Theorem \ref{sinter}. In the tree below, we will show how we proceed from graph to supergraph, by adding vertices. We will underline a letter-number combination if the corresponding graph of Appendix 1 has the smallest eigenvalue less than $\beta$. From $C$ to $D$ we use the fact that $C1$, $C3$ to $C6$ of Appendix 1 are not possible.

\vspace{-0.5cm}

\begin{center}
\begin{figure}[h]
\includegraphics[width=0.9\textwidth]{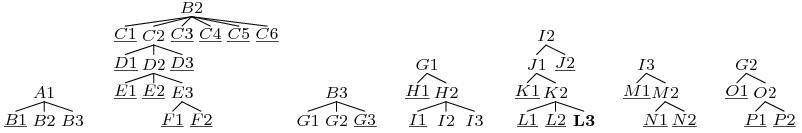}
\end{figure}
\end{center}

\vspace{-0.5cm}

So, this shows this theorem is true under the assumption that it contains a $K_{2,1,1}$ as an induced subgraph.

From now on, we will consider only connected $3$-regular graphs that do not contain $K_{2,1,1}$. From $R1$ to $S$ of Appendix 1, we use the fact that $A1$ and $Q2$ of Appendix 1 do not occur. From there on, we assume that $G$ can not contain $Q2$ as an induced subgraph.

\vspace{-0.5cm}

\begin{center}
\begin{figure}[h]
\includegraphics[width=0.9\textwidth]{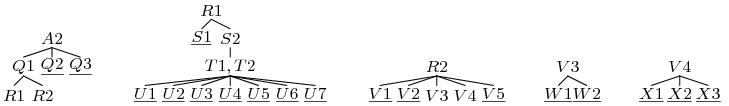}
\end{figure}
\end{center}

\vspace{-0.5cm}

\epf

\begin{remark}
$T5$ and $T7$ are isomorphic.
\end{remark}

\section{Remarks and questions}

\begin{remark}
\begin{itemize}
\item[(i)] Consider $K_{2,3}$ as the $(2,3)$-semi-regular graph in Remark \ref{semi}. Then the graph $G_3$ has $-1-\sqrt2$ as the smallest eigenvalue. This implies that there is a $k$-regular graph with smallest eigenvalue $-1-\sqrt2$ for all $k\geq 3$.

\item[(ii)] Consider three partitions of $\{1, 2, \ldots , 12\}$:
\begin{center}
$P=\{p_i \mid i=1, 2, \ldots , 6\}$ where $p_i=\{i, i+6\}$; \\
$Q=\{q_j \mid j=1, 2, 3\}$ where $q_j=\{j, j+3, j+6, j+9\}$; \\
$R=\{r_l \mid l=1, 2, 3, 4\}$ where $r_l=\{3l-2, 3l-1, 3l\}.$
\end{center}
Then the smallest eigenvalue of the graph which is obtained by a similar argument of construction of $G_k$ in Section \ref{2} is $\alpha_1$. This means that there is a $k$-regular graph with smallest eigenvalue $\alpha_1$ for all $k\geq 4$.
\end{itemize}
\end{remark}

\begin{question}
\begin{itemize}
\item[(i)] It is not known that whether $\hat{\theta_k}$ and $\hat{\eta_k}$ are monotone decreasing sequences or not.
\item[(ii)] It is not known that whether there are infinitely many regular graph with $\lambda_{min}\in [-1-\sqrt2 ,-2)$ which are not $\{\mathfrak{H}_2, \mathfrak{H}_5, \mathfrak{H}_7, \mathfrak{H}_9\}$-line graph.
\end{itemize}
\end{question}

\section*{Acknowledgements}
Part of this work was done while visiting the Graduate School of Information Sciences(GSIS) at Tohoku University. The author greatly appreciates the hospitality of Profs. Munemasa and Obata. And also, the author would like to thank Prof. Koolen and Jongyook Park for the careful reading they did.

\section*{Appendix 1}

\begin{center}
\begin{figure}[h]
\includegraphics[width=0.9\textwidth]{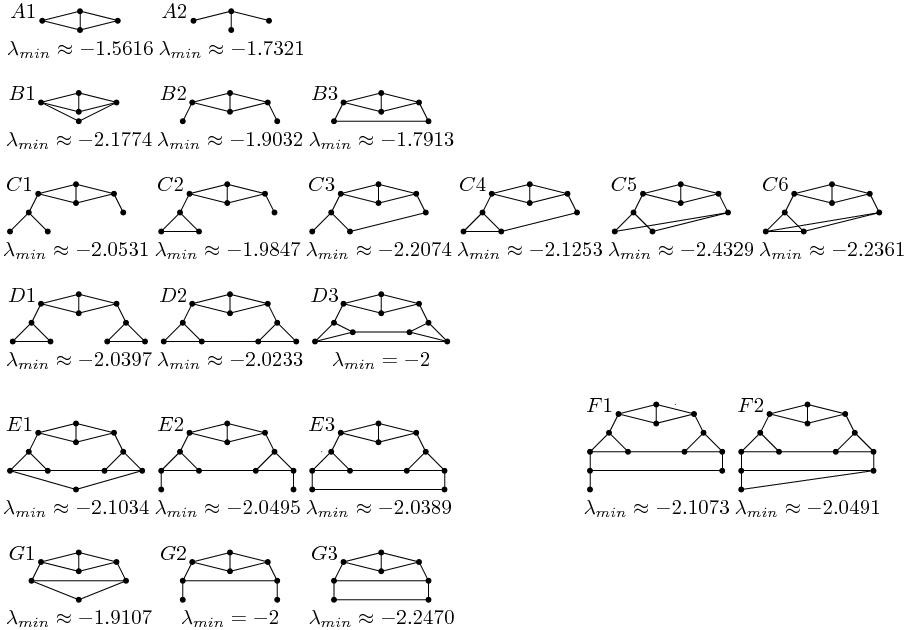}
\end{figure}
\end{center}

\begin{center}
\begin{figure}[h]
\includegraphics[width=0.9\textwidth]{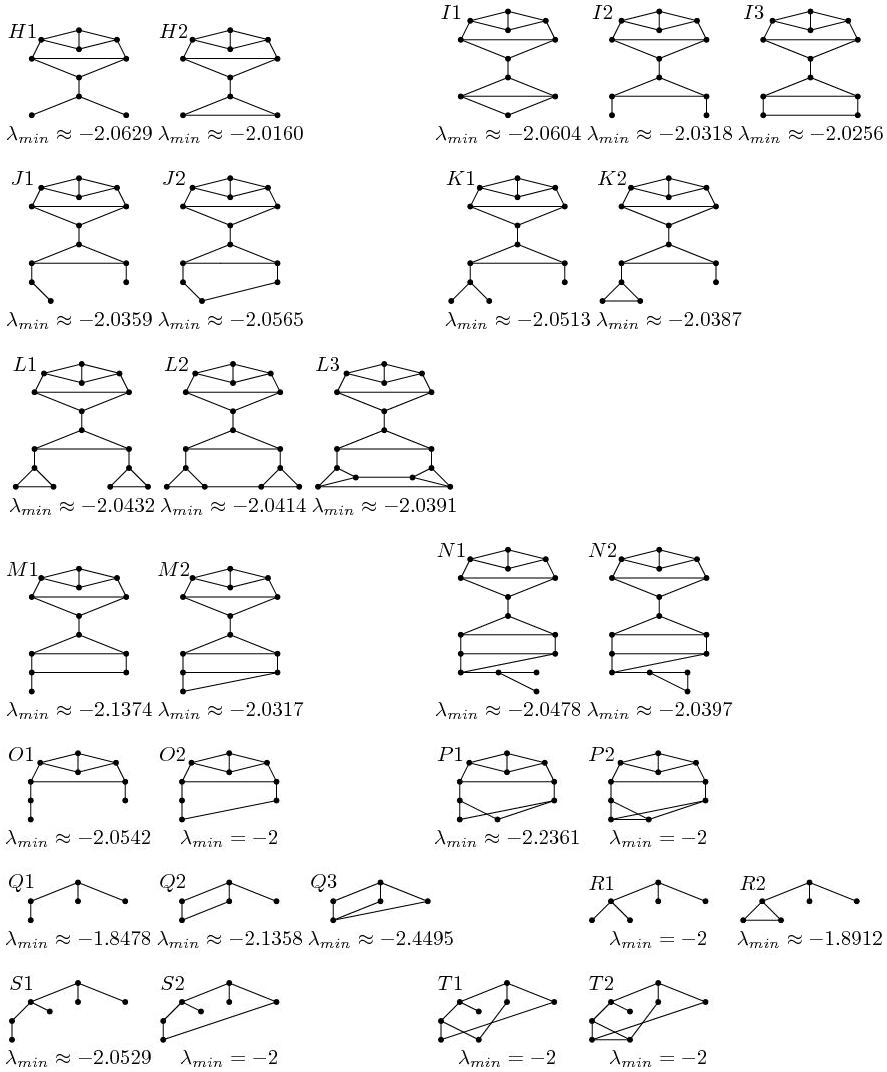}
\end{figure}
\end{center}

\begin{center}
\begin{figure}[h]
\includegraphics[width=0.9\textwidth]{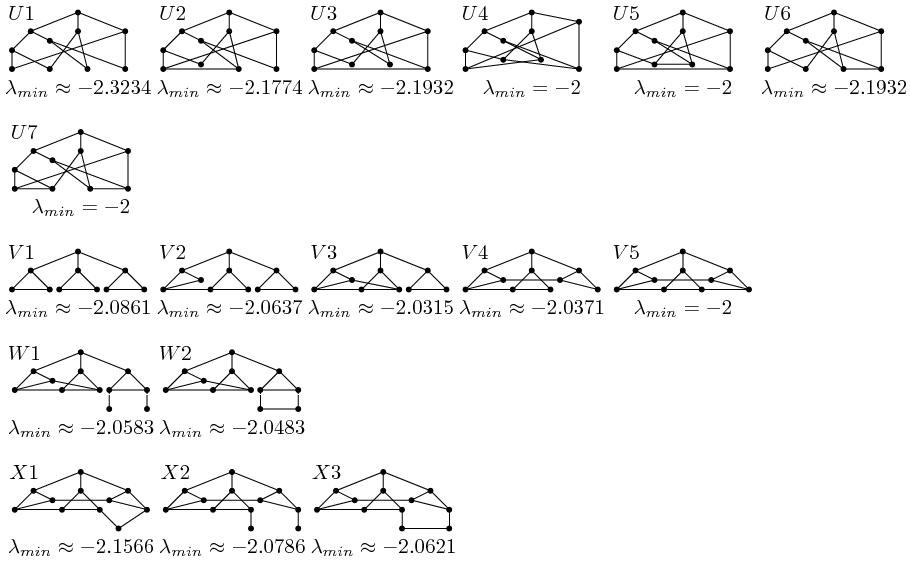}
\end{figure}
\end{center}

\end{document}